\documentclass[a4paper,10pt]{article}
\usepackage{amssymb,amsmath}
\newcommand{\Q}{\mathbb{Q}}
\newcommand{\R}{\mathbb{R}}

\newcommand{\N}{\mathbb{N}}
\newtheorem{theorem}{Theorem}
\newtheorem{lemma}{Lemma}

\newcommand{\Z}{\mathbb{Z}}

\newcommand{\twosum}[2]{\sum_{\substack{#1\\#2}}}

\newcommand{\ep}{\varepsilon}

\begin{document}
\title{Bounds for the Cubic Weyl sum}
\author{D.R. Heath-Brown\\Mathematical Institute, Oxford}
\date{}
\maketitle

\section{Introduction}
In this paper we shall consider bounds for the cubic Weyl sum
\[S(\alpha,N)=\sum_{n\le N} e(\alpha n^3),\]
where $e(x)=\exp(2\pi i x)$ as usual. The classical bound, due
essentially to Weyl \cite{Weyl}, shows that
$S(\alpha,N)\ll_{\ep}N^{3/4+\ep}$ for any $\ep>0$, providing that there
is a rational number $a/q$ with denominator in the range $N\le q\le N^2$,
for which we have $|\alpha-a/q|\le q^{-2}$. It is clear that a condition on
rational approximations to $\alpha$ will be necessary, and the exact
condition here is unimportant. Of greater significance is the
exponent $3/4$ in the Weyl estimate, which has never been improved
on.  An alternative method to bound $S(\alpha,N)$ has been given by
Vaughan \cite[Theorem 3]{Vaughan}, leading to exactly the same exponent
$3/4$.  If $\alpha$ is a real algebraic irrational then Roth's Theorem
shows that the Diophantine approximation condition is met, so that
\[S(\alpha,N)\ll_{\ep,\alpha}N^{3/4+\ep}\]
for all $N$.

The goal of this paper is to show how an improvement can be made for
special values of $\alpha$.  Unfortunately our result depends on
unproved hypothesis, namely the $abc$-conjecture.  This states that if
$\ep>0$ is given there is a constant $K(\ep)$ such that
\[\max\{|a|,|b|,|c|\}\le K(\ep)\left(\prod_{p|abc}p\right)^{1+\ep}\]
for any coprime positive integers $a,b,c$ with $a+b=c$.

We shall then prove the
following bound.
\begin{theorem}\label{sqrt2}
Let $\alpha\in\R-\Q$ be a quadratic irrational.  Assume the truth of the
$abc$-conjecture.  Then
\[S(\alpha,N)\ll_{\ep,\alpha}N^{5/7+\ep}\]
for any $\ep>0$.
\end{theorem}
Note that $5/7=3/4-1/28$.

The underlying idea is to apply the $q$-analogue of van der Corput's
method, which requires a suitable approximation $a/q$ to $\alpha$,
in which $q$ factorizes in a suitable way.  Results of this type were
proved in an Oxford DPhil thesis by Ringrose \cite{Ringrose} in 1985,
but not otherwise published.  We therefore establish a variant of
Ringrose's result here.
\begin{theorem}\label{genthm}
Suppose that $a$ and $q$ are coprime integers with $N\le q\le
N^{3/2}$.  Suppose further that $q=q_1q_2q_3$ with the factors
$q_1,q_2,q_3$ coprime in pairs and $q_3$ square-free.  Then if
$N\le\min\{q_1q_3,q_2q_3\}$ we have
\[S(\alpha,N)\ll_{\ep}
\left(1+N^3\left|\alpha-\frac{a}{q}\right|\right)
(N^{1/2}q_1^{1/2}+N^{1/4}q^{1/4}q_2^{1/4}+
N^{1/4}q^{1/4}q_3^{1/8})q^{\ep},
\]
for any $\ep>0$.
\end{theorem}

Theorem \ref{sqrt2} is an easy consequence of Theorem \ref{genthm},
along with the following result on Diophantine approximation with
smooth denominators.

\begin{theorem}\label{smden}
Let $\alpha\in\R$ be a quadratic irrational, and let $\ep>0$ be
given.  Then there is a constant $C(\alpha,\ep)$ such that, for any
$N\in\N$, one can solve
\[\left|\alpha-\frac{a}{q}\right|\le\frac{C(\alpha,\ep)}{qN},\;\;\;
(a\in\Z,\;\;q\in\N,\;\;q\le N)\]
with $q$ having no prime factors $p>q^{\ep}$.
\end{theorem}
This result provides approximations almost as strong as Dirichlet's
Theorem yields.  However the hypothesis that $\alpha$ is quadratic
makes the proof rather simple. One can quite easily improve the statement
of the theorem slightly to say that any prime power factor $p^e$ of
$q$ has $p^e\le q^{\ep}$.  However we do not need this for our
application.  Unfortunately we are unable to produce values of $q$
which are ``nearly square-free'', and it is at this point that we must
call on the $abc$-conjecture. Indeed the reader may notice that it
suffices for Theorem \ref{sqrt2} that the largest power-full divisor
of $q$ should be at most $q^{11/21}$, but we are unable to produce
numbers $q$ of this type.
\bigskip

This paper was written while the author was attending the trimestre on
Diophantine Equations at the Hausdorff Institute of Mathematics in
Bonn.  The hospitality and financial support of the institute is
gratefully acknowledged.  The idea of using factorization properties
of the sequence $q_n$ to facilitate the $q$-analogue of van der
Corput's method arose in a somewhat different context in a
conversation with Jimi Truelsen.  His contribution is also gratefully
acknowledged.

\section{Preliminary Steps}
If $\delta =\alpha-a/q$ we have
\[S(\alpha,N)=S(\frac{a}{q},N)e(\delta N^3)-\int_0^N(2\pi i\delta)3t^2
S(\frac{a}{q},t)dt\ll (1+N^3|\delta|)\max_{t\le
N}|S(\frac{a}{q},t)|,\]
by partial summation.  Moreover, if we set
\[S(a,h;q)=\sum_{n=1}^qe((an^3+hn)/q)\]
and
\[T(h,t;q)=\sum_{n\le t}e(-hn/q)\]
we have
\[S(\frac{a}{q},t)=q^{-1}\sum_{-q/2<h\le q/2}S(a,h;q)T(h,t;q).\]
Since $S(a,0;q)\ll q^{2/3}$ 
the term $h=0$ contributes $\ll Nq^{-1/3}$, which is satisfactory.

For the remaining terms we note that 
\[T(h,t;q)\ll \min(N,q/|h|),\;\;\;\mbox{and}\;\;\;
\frac{\partial T(h,t;q)}{\partial h}
\ll q^{-1}N\min(N,q/|h|)\]
for $|h|\le q/2$ and $t\le N$. From now on we shall assume that 
\[\left|\sum_{-q/2<h<0}S(a,h;q)T(h,t;q)\right|\le
\left|\sum_{0< h\le q/2}S(a,h;q)T(h,t;q)\right|,\]
the alternative case being treated in exactly the same way.

We proceed to define $K=[q/N]$ and
\[\eta(r)=\max_{0\le L\le K}
\left|\sum_{(r-1)K< h \le (r-1)K+L}S(a,h;q)\right|.\]
We then find by partial summation that
\[\sum_{(r-1)K\le h\le (r-1)K+K'}S(a,h;q)T(h,t;q)\ll 
\eta(r)\min(N,\frac{q}{(r-1)K})\ll \frac{N}{r}\eta(r)\]
for any integer $r\ge 0$ and any $K'\le K$. 
Summing for $r\le q$ (which is more than adequately large) we deduce that
\[\sum_{-q/2<h\le q/2}S(a,h;q)T(h,t;q)\ll N\sum_{r\le q}\frac{\eta(r)}{r}.\] 
We may therefore conclude as follows.
\begin{lemma}\label{first}
With the definitions above, if $N\le q\le N^{3/2}$ we have
\[S(\frac{a}{q},N)\ll \frac{N}{q}\sum_{r\le q}\frac{\eta(r)}{r}.\] 
\end{lemma}
It follows as a special case of Loxton and Vaughan \cite[Theorem 1]{LV} that
\begin{equation}\label{Sb}
S(a,h;q)\ll_{\ep}q^{1/2+\ep}(q,h)^{1/4}
\end{equation}
for any $\ep>0$, whence
\begin{eqnarray*}
\sum_{r\le q}\frac{\eta(r)}{r}&\ll&\sum_{r\le q}\sum_{n\le K}
\frac{|S(a,(r-1)K+n;q)|}{r}\\
&\ll& \sum_{m\le qK}|S(a,m;q)|\min(1,\frac{K}{m})\\
&\ll_{\ep}&q^{1/2+\ep}\sum_{m\le qK}(q,m)^{1/4}\min(1,\frac{K}{m}).
\end{eqnarray*}
We now use the following easy lemma, which we shall prove at the end
of this section.
\begin{lemma}\label{hcf}
Let positive integers $v$ and $H_1\le H_2$ be given.  Then for any 
fixed $\ep>0$ we have
\[\sum_{H_2-H_1<h\le H_2}(h,v)^{\rho}\ll_{\ep}\{H_1+\min(v,H_2)\}v^{\ep}\]
uniformly for $\rho\le 1$, and in particular
\[\sum_{1\le h\le H_2}(h,v)^{\rho}\ll_{\ep}H_2v^{\ep}.\]
\end{lemma}
>From this it follows by partial summation that
\begin{equation}\label{ds}
\sum_{r\le q}\frac{\eta(r)}{r}\ll_{\ep}
q^{1/2+\ep}\sum_{m\le qK}(q,m)^{1/4}\min(1,\frac{K}{m})
\ll_{\ep}Kq^{1/2+2\ep},
\end{equation}
whence
\[S(\frac{a}{q},N)\ll_{\ep}\frac{N}{q}Kq^{1/2+2\ep}\ll_{\ep}q^{1/2+2\ep}
\ll N^{3/4+3\ep},\]
since we are assuming that $q\le N^{3/2}$.  We thus recover the
classical exponent $3/4$.  The argument above is equivalent to that
given by Vaughan, mentioned in the introduction.

To prove Lemma \ref{hcf} we merely note that
\begin{eqnarray*}
\sum_{H_2-H_1<h\le H_2}(h,v)^{\rho}&\le&\sum_{H_2-H_1<h\le H_2}(h,v)\\
&\le&\sum_{d|v,\,d\le H_2}d\#\{H_2-H_1<h\le H_2:\, d|h\}\\
&\le &\sum_{d|v,\,d\le H_2}d\{\frac{H_1}{d}+1\}\\
&\le&H\sum_{d|v}\{H_1+\min(v,H_2)\}\\
&\ll_{\ep}&\{H_1+\min(v,H_2)\}v^{\ep}.
\end{eqnarray*}

\section{The First Iteration}
In order to improve on the classical bound 
we must demonstrate some cancellation amongst
the term $S(a,h;q)$ in the sum $\eta(r)$.  We begin by noting the
factorization property
\begin{equation}\label{prod}
S(a,h;uv)=S(av^2,h;u)S(au^2,h;v)
\end{equation}
for coprime $u$ and $v$.  We begin the van der Corput argument by
writing $\eta(r)=|\Sigma|$, where
\[\Sigma=\sum_{h\in I}S(a,h;q)\]
for an appropriate interval $I\subseteq((r-1)K,rK]$.  Then if
we impose the condition
\begin{equation}\label{cond1}
q_2q_3\ge N,
\end{equation}
we will have $K\ge q_1$.  We then set
$M=[K/q_1]\ge 1$ and observe that
\begin{eqnarray*}
M\Sigma&=&\sum_{m\le M}\sum_{h:\, h+mq_1\in I}S(a,h+mq_1;q)\\
&=&\sum_{(r-2)K<h<rK}\sum_{m\le M:\, h+mq_1\in I}S(a,h+mq_1;q).
\end{eqnarray*}
We now apply (\ref{prod}) with $u=q_1$ and $v=q_2q_3$, so that
\begin{eqnarray*}
S(a,h+mq_1;q)&=&S(a',h+mq_1;q_1)S(b,h+mq_1;q_2q_3)\\
&=&S(a',h;q_1)S(b,h+mq_1;q_2q_3),
\end{eqnarray*}
where we have written $a'=aq_2^2q_3^2$ and $b=aq_1^2$.  It follows
that
\[M\Sigma=\sum_{(r-2)K<h<rK}S(a',h;q_1)\sum_{m\le M:\, h+mq_1\in I}
S(b,h+mq_1;q_2q_3).\]
We now apply Cauchy's inequality to produce
\begin{equation}\label{etab}
M^2\eta(r)^2=M^2|\Sigma|^2\le\eta_1(r)\eta_2(r),
\end{equation}
where
\[\eta_1(r)=\sum_{(r-2)K<h<rK}|S(a',h;q_1)|^2\]
and
\[\eta_2(r)=\sum_{(r-2)K<h<rK}\left|\sum_{m\le M:\, h+mq_1\in I}
S(b,h+mq_1;q_2q_3)\right|^2.\]

We may estimate $\eta_1$ using Lemma \ref{hcf}.  By
(\ref{Sb}) we have
\begin{eqnarray}\label{Sig1}
\eta_1(r)&=&\sum_{(r-2)K<h<rK}|S(a',h;q_1)|^2\nonumber\\
&\ll_{\ep}&q_1^{1+2\ep}\sum_{(r-2)K<h<rK}(h,q_1)^{1/2}\nonumber\\
&\ll_{\ep}&q_1^{1+3\ep}(K+q_1)\nonumber\\
&\ll_{\ep}&q_1^{1+3\ep}K,
\end{eqnarray}
since $K\ge q_1$, as noted above.

To handle $\eta_2(r)$ we expand the square to produce
\begin{eqnarray*}
\eta_2(r)&=&\sum_{m_1,m_2\le M}\twosum{h:\, h+m_1q_1\in I}{h+m_2q_1\in
I}S(b,h+m_1q_1;q_2q_3)\overline{S(b,h+m_2q_1;q_2q_3)}\\
&=&\twosum{n_1,n_2\in I}{q_1|n_1-n_2}
S(b,n_1;q_2q_3)\overline{S(b,n_2;q_2q_3)}N(b_1,b_2),
\end{eqnarray*}
where $N(b_1,b_2)$ is the number of triples
$(h,m_1,m_2)\in\Z\times\N\times\N$ for which $m_1,m_2\le M$ and
$h+m_1q_1=n_1,\,h+m_2q_1=n_2$.  Then 
\[N(b_1,b_2)=M-q^{-1}|n_1-n_2|,\]
whence
\begin{eqnarray*}
\eta_2(r)
&=&\twosum{n_1,n_2\in I}{q_1|n_1-n_2}(M-q^{-1}|n_1-n_2|)
S(b,n_1;q_2q_3)\overline{S(b,n_2;q_2q_3)}\\
&=&\sum_{|m|\le M}(M-|m|)\sum_{n, n+mq_1\in I}
S(b,n+mq_1;q_2q_3)\overline{S(b,n;q_2q_3)}\\
&=&\sum_{|m|\le M}(M-|m|)\sum_{n\in I(m)}S_2(b,m,n;q_2q_3),
\end{eqnarray*}
where $I(m)$ is a subinterval of $((r-1)K,rK]$ given by
\[I(m)=\{x\in\R:\, x,x+mq_1\in I\}\]
and $S_2(b,n;q_2q_3,)$ is the exponential sum
\begin{equation}\label{s2def}
S_2(b,m,n;u)=S(b,n+mq_1;u)\overline{S(b,n;u)}.
\end{equation}
We write
\[\eta_3(r,m)=\sum_{n\in I(m)}S_2(b,m,n;q_2q_3),\]
so that our bound becomes
\begin{equation}\label{Sig2}
\eta_2(r)\ll M\sum_{|m|\le M}|\eta_3(r,m)|=M|\eta_3(r,0)|+
\sum_{1\le |m|\le M}|\eta_3(r,m)|.
\end{equation}
Notice that the only dependence of $\eta_3(r,m)$ on $r$ is through the
interval $I(m)$, which our notation has suppressed.

We now combine Lemma \ref{first} with (\ref{etab}), (\ref{Sig1}) and
(\ref{Sig2}) to deduce that
\begin{eqnarray*}
S(\frac{a}{q},N)^2&\ll& 
(\frac{N}{q})^2(\log q)\sum_{r\le q}\frac{\eta(r)^2}{r}\\
&\ll_{\ep}& N^2q^{-2+\ep}M^{-2}q_1^{1+3\ep}K
\sum_{r\le q}\frac{\eta_2(r)}{r}\\
&\ll_{\ep}& N^2q^{-2+\ep}M^{-2}q_1^{1+3\ep}KM\sum_{r\le q}\sum_{|m|\le M}
\frac{\eta_3(r,m)}{r}\\
&\ll_{\ep}& N^2q^{-2+\ep}q_1^{2+3\ep}\sum_{r\le q}\sum_{|m|\le M}
\frac{|\eta_3(r,m)|}{r},
\end{eqnarray*}
whence
\[S(\frac{a}{q},N)^2\ll_{\ep}T_1+T_2,\]
where
\[T_1=N^2q^{4\ep}(q_2q_3)^{-2}\sum_{r\le q}\frac{|\eta_3(r,0)|}{r}\]
and
\begin{equation}\label{T2def}
T_2=N^2q^{4\ep}(q_2q_3)^{-2}\sum_{r\le q}\;
\sum_{1\le |m|\le M}\frac{|\eta_3(r,m)|}{r}.
\end{equation}
However
\begin{eqnarray*}
\eta_3(r,0)&=&\sum_{n\in I(m)}S_2(b,0,n;q_2q_3)\\
&=&\sum_{n\in I(m)}|S(b,n;q_2q_3)|^2\\
&\ll_{\ep}&(q_2q_3)^{1+2\ep}\sum_{(r-1)K<n\le rK}(q_2q_3,n)^{1/2}
\end{eqnarray*}
by (\ref{Sb}), whence 
\begin{eqnarray*}
T_1&\ll_{\ep}&N^2q^{4\ep}(q_2q_3)^{-2}\sum_{r\le q}r^{-1}
(q_2q_3)^{1+2\ep}\sum_{(r-1)K<n\le rK}(q_2q_3,n)^{1/2}\\
&\ll_{\ep}&N^2q^{6\ep}(q_2q_3)^{-1}\sum_{n\le qK}\frac{K}{n}
(q_2q_3,n)^{1/2}\\
&\ll_{\ep}&N^2q^{6\ep}(q_2q_3)^{-1}Kq^{\ep}\\
&\ll_{\ep}&Nq^{7\ep}q_1\\
\end{eqnarray*}
using Lemma \ref{hcf}.  This is satisfactory for Theorem \ref{genthm}.

We summarize the state of play as follows.
\begin{lemma}\label{1stit}
When $q_2q_3\le N$ and $q\le N^{3/2}$ we have
\[S(\frac{a}{q},N)^2\ll_{\ep}Nq^{7\ep}q_1+T_2\]
with $T_2$ as in {\rm (\ref{T2def})}.
\end{lemma}

This completes the first application of the van der Corput
``A-process''.  We can check that nothing of significance has been
lost at this stage.  Thus if we use the bound (\ref{Sb}), ignoring the
highest common factor terms for simplicity, we would get a bound
\[\ll_{\ep}M\sum_{|m|\le M}\sum_{n\in I(m)}(q_2q_3)^{1+2\ep}
\ll_{\ep}K^3q_1^{-2}(q_2q_3)^{1+2\ep}\]
for $\eta_2(r)$.  Combining this with (\ref{etab}) and (\ref{Sig1})
would lead to $\eta(r)\ll_{\ep}Kq^{(1+3\ep)/2}$, allowing us to recover
the bound (\ref{ds}).  As previously observed this in turn would lead
to the classical exponent $3/4$ for the original sum $S(a/q,N)$.

Although nothing has been lost, our manipulations have produced an
advantage. We need to demonstrate that there is some cancellation in
the sum $\eta_3(r,m)$. The range for $n$ in this sum
is of the same kind as in the previous sum $\Sigma$,
and the exponential sum $S_2(b,m,n;q_2q_3)$ which occurs is more
complicated than before, but crucially the modulus $q_2q_3$ for the
exponential sum is smaller than in $\Sigma$, where it was $q$.
Unfortunately the modulus is still too large, so that a second
iteration of the A-process is necessary.

\section{The Second Iteration}

For the second iteration we impose the condition
\begin{equation}\label{cond2}
q_1q_3\ge N,
\end{equation}
whence $q_2\le q/N$.  It follows that $q_2\le K$.  We now 
set $U=[K/q_2],$ so that $U\ge 1$. For the sum (\ref{s2def}) the
product formula (\ref{prod}) leads to the relation
\[S_2(b,m,n;q_2q_3)=S_2(b',m,n;q_2)S_2(c,m,n;q_3),\]
where $b'=bq_3^2$ and $c=bq_2^2$. Then, by the arguments leading to
(\ref{etab}) and (\ref{Sig2}), we find that
\begin{equation}\label{Sig3}
|\eta_3(r,m)|^2\le U^{-2}\eta_4(r,m)\eta_5(r,m),
\end{equation}
with
\[\eta_4(r,m)=\sum_{(r-2)K<h<rK}|S_2(b',m,h;q_2)|^2\]
and
\[\eta_5(r,m)=\sum_{|u|\le U}(U-|u|)
\sum_{n\in I(m,u)}S_3(c,m,u,n;q_3).\]
Here $I(m,u)$ is an appropriate subinterval of $((r-1)K,rK]$, and
\[S_3(c,m,u,n;v)=S_2(c,m,n+uq_2;v)\overline{S_2(c,m,n;v)}.\]
By (\ref{Sb}) we will have
\begin{eqnarray*}
S_2(b',m,h;q_2)&\ll_{\ep}&q_2^{1+2\ep}(h+mq_1,q_2)^{1/4}(h,q_2)^{1/4}\\
&\ll_{\ep}&q_2^{1+2\ep}\{(h+mq_1,q_2)^{1/2}+(h,q_2)^{1/2}\}.
\end{eqnarray*}
It follows that
\begin{eqnarray}\label{Sig4}\nonumber
\eta_4(r,m)&\ll_{\ep}&
q_2^{2+4\ep}\sum_{(r-2)K<h<rK}\{(h+mq_1,q_2)+(h,q_2)\}\nonumber\\
&\ll_{\ep}&q_2^{2+4\ep}\sum_{(r-3)K<h<(r+1)K}(h,q_2)\nonumber\\
&\ll_{\ep}&q_2^{2+5\ep}(q_2+K)\nonumber\\
&\ll_{\ep}&q_2^{2+5\ep}K,
\end{eqnarray}
by Lemma \ref{hcf}.  We now set
\[\eta_6(r,m,u)=\sum_{n\in I(m,u)}S_3(c,m,u,n;q_3),\]
whence
\begin{equation}\label{Sig5}
\eta_5(r,m)\ll U\sum_{|u|\le U}|\eta_6(r,m,u)|.
\end{equation}
It now follows that
\begin{eqnarray}\label{T22}
T_2^2&\ll_{\ep}&N^4q^{8\ep}(q_2q_3)^{-4}\left\{M\sum_{r\le
    q}r^{-1}\right\}\left\{\sum_{r\le q}\;
\sum_{1\le |m|\le M}\frac{|\eta_3(r,m)|^2}{r}\right\}\nonumber\\
&\ll_{\ep}&N^4q^{8\ep}(q_2q_3)^{-4}\left\{q_2q_3N^{-1}q^{\ep}\right\}
\times\nonumber\\
&&\hspace{2cm}\mbox{}\times\left\{U^{-1}Kq_2^{2+5\ep}\sum_{r\le q}\;
\sum_{1\le |m|\le M}\;
\sum_{|u|\le U}\frac{|\eta_6(r,m,u)|}{r}\right\}\nonumber\\
&\ll_{\ep}&N^3q^{14\ep}q_3^{-3}\sum_{r\le q}\;
\sum_{1\le |m|\le M}\;\sum_{|u|\le U}\frac{|\eta_6(r,m,u)|}{r}.
\end{eqnarray}
When $u=0$ we have
\begin{eqnarray*}
S_3(c,m,0;q_3)&=&|S_2(c,m.n;q_3)|^2\\
&=&|S(c,n+mq_1;q_3)|^2|S(c,n;q_3)|^2\\
&\ll_{\ep}&q_3^{2+4\ep}(q_3,n+mq_1)^{1/2}(q_3,n)^{1/2}
\end{eqnarray*}
by (\ref{Sb}). The contribution to (\ref{T22}) from terms with $u=0$
is then
\begin{eqnarray*}
&\ll_{\ep}&N^3q^{18\ep}q_3^{-1}\sum_{r\le q}\;
\sum_{1\le |m|\le M}\;\sum_{(r-1)K<n,n+mq_1\le rK}
\frac{(q_3,n+mq_1)^{1/2}(q_3,n)^{1/2}}{r}\\
&\ll_{\ep}&N^3q^{18\ep}q_3^{-1}\sum_{r\le q}\;
\sum_{1\le |m|\le M}\;\sum_{(r-1)K<n,n+mq_1\le rK}
\frac{(q_3,n+mq_1)+(q_3,n)}{n/K}\\
&\ll_{\ep}&N^3q^{18\ep}q_3^{-1}KM\sum_{r\le q}\;
\sum_{(r-1)K<n\le rK}\frac{(q_3,n)}{n}\\
&\ll_{\ep}&N^3q^{18\ep}q_3^{-1}KM
\sum_{n\le qK}\frac{(q_3,n)}{n}\\
&\ll_{\ep}&N^3q^{19\ep}q_3^{-1}KM\\
&\ll_{\ep}&Nq^{19\ep}q_1q_2^2q_3,
\end{eqnarray*}
by Lemma \ref{hcf}

We complete the second van der Corput A-process by combining this with
Lemma \ref{1stit} and (\ref{T22}) to deduce the following bound.
\begin{lemma}\label{2ndit}
When $q_2q_3, q_1q_3\ge N$ and $q\le N^{3/2}$ we have
\[S(\frac{a}{q},N)^4\ll_{\ep}q^{19\ep}\{N^2q_1^2+Nq_1q_2^2q_3
+N^3q_3^{-3}\sum_{r\le q}\;
\sum_{1\le |m|\le M}\;\sum_{1\le |u|\le U}\frac{|\eta_6(r,m,u)|}{r}\}.\]
\end{lemma}
The first two terms here are suitable for Theorem \ref{genthm}.

\section{The van der Corput B-Process}

To complete the van der Corput argument we will estimate
$\eta_6(r,m,u)$.  We have
\[\eta_6(r,m,n)=q_3^{-1}\sum_{-q_3/2<t\le q_3/2}S_4(c,m,u,t;q_3)
\sum_{n\in I(m,u)}e(-nt/q_3),\]
where
\[S_4(c,m,u,t;v)=\sum_{n=1}^vS_3(c,m,u,n;v)e(nt/v).\]
Since $I(m,n)$ is an interval of length at most $K$ this leads to
\begin{equation}\label{e6e}
\eta_6(r,m,u)\ll q_3^{-1}\sum_{-q_3/2<t\le q_3/2}
\min(K\,,\,\frac{q_3}{|t|})|S_4(c,m,u,t;q_3)|.
\end{equation}
The sum $S_4(c,m,u,n,t;v)$ has a multiplicative property
\begin{equation}\label{m4}
S_4(c,m,u,t;vw)=S_4(cw^2,m,u,\overline{w}t;v)S_4(cv^2,m,u,\overline{v}t;w),
\end{equation}
where $w\overline{w}\equiv 1\mod{v}$ and $v\overline{v}\equiv 1\mod{w}$.
It therefore suffices to bound sums to prime-power modulus.  Indeed,
since we are assuming $q_3$ to be square-free it will be enough to
consider the case in which the modulus is prime.  It would be good if
we were able to remove the square-freeness condition, by handling
$S_4$ for prime power moduli, but this appears to be unduly complicated.

In view of the definition of $S_3(c,m,u,n;v)$ we find that
\begin{eqnarray*}
S_4(c,m,u,t;v)&=&\sum_{n=1}^v S^{(1)}\overline{S^{(2)}S^{(3)}}S^{(4)}e(nt/v)\\
&=&\sum_{w,x,y,z=1}^v\sum_{n=1}^ve(f(w,x,y,z,n)/v),
\end{eqnarray*}
where
\[S^{(1)}=S(c,n+uq_2+mq_1;v),\;\;\;S^{(2)}=S(c,n+uq_2;v),\]
\[S^{(3)}=S(c,n+mq_1;v),\;\;\;S^{(4)}=S(c,n;v)\]
and
\begin{eqnarray*}
f(w,x,y,z,n)&=&c(w^3-x^3-y^3+z^3)+w(n+uq_2+mq_1)\\
&&\hspace{1cm}\mbox{}-x(n+uq_2)-y(n+mq_1)+zn+tn.
\end{eqnarray*}
When we perform the summation over $n$ this produces
\[v\twosum{w,x,y,z\!\!\!\!\!\pmod{v}}{v|w-x-y+z+t}
e\{\frac{c(w^3-x^3-y^3+z^3)+uq_2(w-x)+mq_1(w-y)}{v}\}.\]
We substitute $z=x+y-w-t$ so that the denominator becomes
\begin{eqnarray*}
\lefteqn{c(w^3-x^3-y^3+(x+y-w-t)^3)+uq_2(w-x)+mq_1(w-y)}\\
&=& 3c(x+y)(w-x)(w-y)-3ct(x+y-w)^2+3ct^2(x+y-w)-ct^3\\
&&\hspace{1cm}\mbox{}+uqw_2(w-x)+mq_1(w-y)\\
&=&3cWXY-\frac{3}{4}ct(W+X+Y)^2+\frac{3}{2}ct^2(W+X+Y)-ct^3\\
&&\hspace{1cm}\mbox{}-uq_2X-mq_1Y,
\end{eqnarray*}
on writing $W=x+y$ and $X=x-w,\,Y=y-w$.
Thus if $v$ is a prime $p$
which does not divide $6t$ we find that
\[=p\sum_{W,X,Y\!\!\!!\!\pmod{p}}e(g(W,X,Y)/p)\]
where $g(W,X,Y)$ takes the shape
\[c'WXY+t'(W^2+W^2+Y^2+2XY+2WX+2WY)+\mu_1 W+\mu_2 X+\mu_3 Y-ct^3\]
with $p\nmid c't'$.
Exponential sums of this type have been treated by Bombieri and
Sperber \cite[Theorem 7]{BS}.  Their work shows that
\[S_4(c,m,u,t;p)\ll p^{5/2}\]
for such primes.  When $p\mid t$ it is easy to see that $S_4(c,m,u,t;p)\ll
p^2(p,m,u)$, whence in general we have $S_4(c,m,u,t;p)\ll
p^{5/2}(p,t,m,u)^{1/2}$ for all primes $p$.
While it is convenient to call on a theorem from the
literature, we remark that it is possible to evaluate $S_4(c,m,u,t;p)$
in terms of exponential sums in one variable, for which it suffices to
use Weil's theorem rather than Deligne's.

By the multiplicative property (\ref{m4}) we now deduce that
\[S_4(c,m,u,t;q_3)\ll_{\ep}q_3^{5/2+\ep}(q_3,t,m,u)^{1/2},\]
whence (\ref{e6e}) and Lemma \ref{hcf} yield
\begin{eqnarray*}
\eta_6(r,m,u)&\ll_{\ep}& q_3^{3/2+\ep}\sum_{-q_3/2<t\le q_3/2}
\min(K,\frac{q_3}{|t|})(q_3,t,m,u)^{1/2}\\
&\ll_{\ep}& q_3^{3/2+\ep}\{K(q_3,m,u)^{1/2}+
\sum_{1\le t\le q_3/2}\frac{q_3}{|t|}(q_3,t)^{1/2}\}\\
&\ll_{\ep}& q_3^{3/2+\ep}\{K(q_3,m,u)^{1/2}+q_3^{1+\ep}\},
\end{eqnarray*}
by partial summation.  It then follows using Lemma \ref{hcf} that
\begin{eqnarray*}
\lefteqn{N^3q_3^{-3}\sum_{r\le q}\;
\sum_{1\le |m|\le M}\;\sum_{1\le |u|\le U}\frac{|\eta_6(r,m,u)|}{r}}
\hspace{2cm}\\
&\ll_{\ep}&N^3q_3^{-3/2+\ep}\{KMUq_3^{\ep}+MUq_3^{1+2\ep}\}\\
&\ll_{\ep}& q^{3\ep}\{q_1^2q_2^2q_3^{3/2}+Nq_1q_2q_3^{3/2}\}\\
&\ll_{\ep}& Nq_1q_2q_3^{3/2}q^{3\ep},
\end{eqnarray*}
since
\[N=\left(\frac{N^{3/2}}{q_1q_2}\right)^2\frac{q_1^2q_2^2}{N^2}\ge
\left(\frac{q}{q_1q_2}\right)^2\frac{q_1^2q_2^2}{N^2}=
q_1q_2\frac{q_1q_3}{N}\cdot\frac{q_2q_3}{N}\ge q_1q_2.\]
Now, when we insert this estimate into Lemma \ref{2ndit}, we see that
Theorem \ref{genthm} follows, with a new value for $\ep$.

\section{Proof of Theorems \ref{smden} and   \ref{sqrt2}} 
To prove Theorem \ref{smden} it clearly suffices to suppose that
$\alpha=\sqrt{d}$ for some non-square $d\in\N$.  Let $a,b\in\N$ be
solutions of the Pell equation $a^2-db^2=1$ and define $p_n,q_n$ by
\[p_n+q_n\sqrt{d}=\eta^n\]
where $\eta=a+b\sqrt{d}$.  It follows that
\[q_n=(2\sqrt{d})^{-1}(\eta^n-\eta^{-n})=
(2\sqrt{d})^{-1}\prod_{k|n}\Phi_k(\eta,\eta^{-1}),\]
where
\[\Phi_k(X,Y)=\prod_{1\le h\le k,\,(h,k)=1}(X-e^{2\pi ih/k}Y)\]
is the $k$-th cyclotomic polynomial.  Thus $\Phi_k(X,Y)\in\Z[X,Y]$ and
$\Phi_k(X,X^{-1})=\Phi_k(X^{-1},X)$ except for $k=1$, in which case
$\Phi_1(X,Y)=X-Y$. It follows that
\[q_n=b\prod_{k|n,\, k\ge 2}r_k\]
for integers $r_k=|\Phi_k(\eta,\eta^{-1})|$.  Moreover
\[r_k\le(\eta+\eta^{-1})^{\phi(k)}=(2a)^{\phi(k)}\le (2a)^{\phi(n)}.\]

Now fix an integer $m$ such that
\[\frac{\phi(m)}{m}\le\ep\frac{\log \eta}{2\log(2a)}.\]
Then if $m|n$ we have $\phi(n)/n\le\phi(m)/m$, whence
\[r_k\le(2a)^{\phi(n)}\le \{(2a)^n\}^{\phi(m)/m}
\le \{(2a)^n\}^{\ep(\log\eta)/(2\log(2a)}=(\eta^{n/2})^{\ep}.\]
However
\[q_n=(2\sqrt{d})^{-1}(\eta^n-\eta^{-n})\ge
(4\sqrt{d})^{-1}\eta^n\ge\eta^{n/2}\]
for large enough $n$, whence $r_k\le q_n^{\ep}$.  It follows that
every prime factor of $q_n$ is at most $q_n^{\ep}$, if $m|n$.

Finally we observe that $q_{n+1}\ll q_n$ with an implied constant
depending on the choice  of $a,b$ and $d$, so that there is a value of
$n$ which is a multiple of $m$ and for which $N\ge q_n\gg N$.  Since
\[|p_n-q_n\sqrt{d}|=\frac{1}{p_n+q_n\sqrt{d}}\le\frac{1}{q_n}\]
the theorem then follows.
\bigskip

To deduce Theorem \ref{sqrt2} we write $\alpha=(f+g\sqrt{d})/c$ and
approximate $\sqrt{d}$ as above, with
\[|\sqrt{d}-\frac{u}{v}\le\frac{1}{vV},\;\;\; V\ge v\gg V,\]
where we choose $V=[N^{3/2}/c]$.  Let $v_0$ be the product of all
prime powers $p^e||q$ for which $e\ge 2$. Then the product of the
prime divisors of $v_0$ can be at most $v_0^{1/2}$.  It follows that
\[\prod_{p|u^2v^2d}p\le u\frac{v}{v_0}v_0^{1/2}d.\]
Since $1+v^2d=u^2$ the $abc$-conjecture would
imply that $u\ll_{\ep,\alpha}(uvv_0^{-1/2})^{1+\ep}$, whence
$v_0\ll_{\ep,\alpha}v^{2\ep}$. 

We now write $a_1=fv+gu$ and $q_1=cv$, and set $a=a_1/(a_1,q_1)$ and
$q=q_1/(a_1,q_1)$, so that $a$ and $q$ are coprime, with 
\[q\le q_1=cv\le cV\le N^{3/2}\]
and
\[q\gg_{\alpha} q_1\gg_{\alpha} v\gg_{\ep,\alpha}V\gg_{\alpha} N^{3/2}.\]
Then
\[\left|\alpha-\frac{a}{q}\right|\ll_{\alpha}
\left|\sqrt{d}-\frac{u}{v}\right|\le\frac{1}{vV}\ll_{\ep,\alpha}
\frac{1}{qN^{3/2}}.\]
Moreover every prime factor of $q$ is $O_{\ep,\alpha}(q^{\ep})$, and
if $q_0$ is the product of all $p^e||q$ with $e\ge 2$, then
$q_0\ll_{\ep,\alpha}(q^{\ep})$. 

We proceed to build up coprime square-free divisors $q_2,q_3$ of
$q/q_0$, one prime factor at a time, to produce products in the ranges
\[q^{5/21}\le q_2\ll_{\ep,\alpha}q^{5/21+\ep},\;\;\;
q{10/21}\le q_3\ll_{\ep,\alpha}q^{10/21+\ep}.\]
We will then have $q=q_1q_2q_3$ with
\[q^{2/7-2\ep}\ll_{\ep,\alpha}q_1\le q^{2/7}.\]
One may then verify that the hypotheses of Theorem \ref{genthm} are
satisfied, and that
\[S(\alpha,N)\ll_{\ep,\alpha}q^{10/21+2\ep}\ll_{\ep,\alpha}N^{5/7+3\ep}.\]
This suffices for Theorem \ref{sqrt2}, on re-defining $\ep$.

\bigskip
\bigskip

Mathematical Institute,

24--29, St. Giles',

Oxford

OX1 3LB

UK
\bigskip

{\tt rhb@maths.ox.ac.uk}

\end{document}